\documentclass[smallextended]{article}
\usepackage{amsmath,amssymb,amsthm}

\title{On a relation between certain character values of symmetric groups and its connection with creation operators of symmetric functions}
\author{
Masaki Watanabe \\ 
Graduate School of Mathematical Sciences, The University of Tokyo, \\
3-8-1 Komaba Meguro-ku Tokyo 153-8914, Japan \\ \texttt{mwata@ms.u-tokyo.ac.jp}}
\date{\empty}

\newcommand{\kap}{\kappa}
\newcommand{\lmb}{\lambda}
\newcommand{\Lmb}{\Lambda}
\newcommand{\emp}{\varnothing}
\newcommand{\surj}{\twoheadrightarrow}

\newcommand{\mor}{\rightarrow}
\newcommand{\QQ}{\mathbb{Q}}
\newcommand{\ZZ}{\mathbb{Z}}
\newcommand{\CC}{\mathbb{C}}

\newcommand{\QB}{\mathcal{B}}

\newtheorem{lem}{Lemma}[section]
\newtheorem{cor}{Corollary}[section]
\newtheorem{prop}{Proposition}[section]
\newtheorem{defn}{Definition}[section]
\newtheorem{thm}{Theorem}[section]
\newtheorem{tmplem}{Lemma}[section]

\theoremstyle{definition}
\newtheorem{rem}{Remark}

\newtheorem{eg}{Example}

\newtheorem{prf}{Proof}

\begin{document}
\maketitle

In this paper, we derive a relation of new kind between certain character values of symmetric groups in terms of so-called maya diagrams. 
We also investigate a relation between our result and Bernstein's creation operators for Schur functions, 
and consider analogous relations for projective characters of symmetric groups through creation operators for Schur $Q$-functions. 
We also consider analogous relations for characters of Brauer algebras and walled Brauer algebras. 


\section{Introduction}
Irreducible characters of symmetric groups have long been the subject of combinatorial representation theory. 
Still one sometimes encounters curious phenomena. 
In examining the character tables of symmetric groups, we noticed a curious correspondence 
between the irreducible character values of $S_n$ at transpositions 
and the degrees of the irreducible characters of $S_{n-2}$ for $n \leq 7$, which is exhibited in \S2.2. 
For example, the degrees of irreducible characters of $S_4$ are $1,3,2,3,1$, 
while the values of irreducible characters of $S_6$ at transpositions are $1,3,2,3,1,0,-1,-3,-2,-3,-1$. 
Even though it does not continue beyond $n \leq 7$, this looks too nice to be a sheer coincidence. 
One of the purposes of this paper is to find a mechanism lying behind this phenomenon, 
which also shows that the phenomenon actually extends beyond such limitations on $n$, 
and to more relationships between entries in other columns of the character tables, 
but in a somewhat less apparent and less straightforward manner. 
Our main result gives an expansion of the irreducible character value $\chi_\lmb(\mu \cup (m))$ 
into a linear combination of the values $\chi_\kappa(\mu)$ ($|\kappa|=|\lmb|-m$) 
for a partition $\mu$ with no parts divisible by $m$. 
In particular, if $m=2$ it can be shown that, for $|\lmb| \leq 7$, we have at most one term in the expansion, 
and setting $\mu=(1, \ldots, 1)$ it gives an explanation for the phenomenon. 
One may notice that the classical Murnaghan-Nakayama formula also matches the description of our result above. 
In fact, we first thought that the phenomenon in question must be explained by a direct application of the Murnaghan-Nakayama formula, 
but, as is explained in \S 2.2, this is not the case. 

We also found that the operators appearing in our formula can be written in a form similar to the creation operators for Schur functions, 
and by considering, instead, the creation operators for Schur $Q$-functions, we arrive at a projective analogue of our result. 

As symmetric groups are in duality with (i.e. constitute the centralizer of) general linear groups, 
Brauer algebras and walled Brauer algebras are also in duality with some classical groups: 
orthogonal and symplectic groups for the former and general linear groups for the latter. 
Thus, they also have ``Frobenius formulae'' for their characters. 
We found that our method also works for these algebras and implies analogous relations between characters as in the symmetric-group case.

The paper is organized as follows. 
In section 2 we review some basic facts about irreducible characters of symmetric groups, 
and derive a relation formula for them, which is our main result. 
In section 3, we relate results in section 2 with Bernstein's creation operators for Schur functions. 
In section 4, we investigate a ``projective analogue'' of the result by considering analogous modification as in section 3 of creation operators for Schur $Q$-functions. 
In section 5 we review facts about Brauer and walled Brauer algebras, 
and see analogous relations as the one derived for symmetric groups in section 2 hold for their irreducible characters. 


\section{Irreducible characters of symmetric groups}
\subsection{Preliminaries}
A sequence of positive integers $\lmb=(\lmb_1, \ldots, \lmb_l)$ with $\lmb_1 \geq \cdots \geq \lmb_l > 0$ is called a \textit{partition}. 
$l$ is called the \textit{length} of $\lmb$ and denoted by $\ell(\lmb)$. The terms $\lmb_i$ are called the \textit{parts} of this partition. 
We write $|\lmb|$ for $\lmb_1+\cdots+\lmb_l$, and if $|\lmb|=n$ then $\lmb$ is called a \textit{partition of $n$} and $n$ is called the \textit{size} of $\lmb$. 
The empty sequence is the unique partition of $0$ (or of length $0$) and is denoted by $\emp$. 
We sometimes write nonempty partitions in multiplicative form: $4^231^3$ stands for the partition $(4,4,3,1,1,1)$, for example. 

Let $S_n$ denote the symmetric group over $n$ letters. It is well known that, over a field of characteristic zero, the isomorphism classes of irreducible linear representations of $S_n$ are indexed by the partitions of $n$. 
In this paper, we use objects called \textit{maya diagrams} rather than partitions to index irreducible representations of symmetric groups. 

\begin{defn}
A maya diagram is a sequence of integers $[x_1, x_2, \ldots]$ with $x_1>x_2>\cdots$ and $x_i = -i$ for $i \gg 1$. 
\end{defn}

We write maya diagrams with $[\,]$ and partitions with $(\,)$ (or with no parentheses if they are in multiplicative form), in order to avoid confusion between these two kinds of objects.

There is an easy correspondence between maya diagrams and partitions, 
by $(\lmb_1, \cdots, \lmb_l) \mapsto [\lmb_1-1, \lmb_2-2, \ldots]$ 
and $[x_1, x_2, \ldots] \mapsto (x_1+1, x_2+2, \cdots)$ 
where partitions are identified with infinite sequences obtained from them by attaching infinitely many zeroes. 
The maya diagram corresponding to a partition $\lmb$ under these bijections is also denoted by $\lmb$, 
but we believe that this causes no confusion. 
The notion of size is also defined for maya diagrams through the correspondence above: 
if $|\lmb|=n$, the corresponding maya diagram $[x_1, x_2, \ldots]$ satisfies $\sum_{i} (x_i+i) = n$. 
We also have $x_i=-i$ for $i > |\lmb|$. 

The following lemma can be easily shown: 
\begin{lem}
For a maya diagram $\lmb=[x_1, x_2, \ldots]$, 
let $y_1<y_2<\cdots$ be the integers not appearing in $x_1, x_2, \ldots$. 
Then the partitions corresponding to $\lmb$ and $[-1-y_1, -1-y_2, \ldots]$ are conjugate to each other. 
In particular, $y_i=i-1$ for $i > |\lmb|$ and $|\lmb|=\sum_{i} (i-1-y_i)$. 
\label{mayalem}
\end{lem}

Let $F$ be a $\CC$-vector space having the set of all maya diagrams as a basis. 

We also define a \textit{charged maya diagram} as a sequence of integers $[x_1, x_2, \ldots]$ with $x_1>x_2>\cdots$ and $x_{i+1}=x_i-1$ for $i \gg 1$. 
Let $\tilde{F}$ be a vector space over $\CC$ having the set of all charged maya diagrams as a basis. 

\begin{rem}
The spaces $F$ and $\tilde{F}$ are sometimes called an infinite or semi-infinite wedge space (\cite{KP, KR}): 
in that case $[x_1, x_2, \ldots]$ is denoted as $v_{x_1} \wedge v_{x_2} \wedge \cdots$. 
\end{rem}

We sometimes use sequences $[x_1, x_2, \ldots]$ which satisfy $x_i=-i$ (or $x_{i+1}=x_i-1$) for $i \gg 1$ but are not necessarily decreasing. 
In such case we consider them as elements of $F$ (or $\tilde{F}$ respectively) by the rule $[\ldots, i, \ldots, j, \ldots]+[\ldots, j, \ldots, i, \ldots]=0$
(in particular, sequences with duplicate terms are equal to zero). 

Let $\chi_{\lmb}$ denote the irreducible character of a symmetric group indexed by a maya diagram or a partition $\lmb$. 

We say an element $w \in S_n$ has \textit{cycle type} $\mu=(\mu_1, \cdots, \mu_l)$ if 
$w$ is a product of disjoint cycles with length $\mu_1 \geq \mu_2 \geq \ldots \geq \mu_l \geq 1$. 
It is well known that two elements in $S_n$ are conjugate if and only if they have the same cycle type. 
Let $\chi_{\lmb}(\mu):=\chi_{\lmb}(w)$ for $w \in S_n$ with cycle type $\mu$. 

Then well-known formulae of Frobenius' and Murnaghan's essentially state the following, under the present setting: 

\begin{thm}
For $r \in \ZZ$, define a $\CC$-linear map $A_r: F \mor F$ by $A_r[x_1, x_2, \ldots, ]=\sum_{i \geq 1}[x_1, x_2, \ldots, x_i-r, \ldots]$ for a maya diagram $[x_1, x_2, \ldots]$
(note that only finitely many summands on the right-hand side are nonzero). 
Then for a maya diagram $\lmb$ and a partition $\mu=(\mu_1, \cdots, \mu_l)$ with $|\lmb|=|\mu|=n$, 
\begin{equation*}
A_{\mu}\lmb := A_{\mu_1} \cdots A_{\mu_l}\lmb = \chi_{\lmb}(\mu)\emp .
\end{equation*}
\label{mnformula}
\end{thm}

\subsection{A motivating phenomenon}
The character tables of symmetric groups can be calculated by Murnaghan's formula explained above. 
For example, it is easy to check that the irreducible characters of $S_4$ have degrees $1, 3, 2, 3, 1$, 
while the irreducible characters of $S_6$ takes values $1, 3, 2, 3, 1, 0, -1, -3, -2, -3, -1$ at transpositions. 
Here one can see an obvious correspondence: 
the degrees of the irreducible characters of $S_4$ as well as their negatives appear again as irreducible character values of $S_6$ at transpositions. 
This correspondence also happens for $S_5$ and $S_7$: the irreducible characters of $S_5$ have degrees $1, 4, 5, 6, 5, 4, 1$ and they appear again as irreducible character values of $S_7$ at transpositions. 
One can also see that it also happens for smaller symmetric groups. 
(Unfortunately, there is no obvious correspondence for $S_6$ and $S_8$: 
for example, $\chi_{4^2}(21^6)=4$ but no irreducible character of $S_6$ has degree $4$, 
and $S_6$ has four irreducible characters of degree $5$ but only two irreducible characters of $S_8$ take value $5$ at transpositions. )

We easily see that Murnaghan's formula does not give, at least directly, the explanation we expect.  
For example, for the character values of $S_6$ on transpositions what Murnaghan's formula gives is as follows: 
\begin{align*}
\chi_{6}(21^4) &=\chi_{4}(1^4)=1,  \\
\chi_{51}(21^4) &=\chi_{31}(1^4)=3,  \\
\chi_{42}(21^4) &=\chi_{4}(1^4)+\chi_{2^2}(1^4)=1+2=3,  \\
\chi_{41^2}(21^4) &=\chi_{21^2}(1^4)-\chi_{4}(1^4)=3-1=2, \text{and} \\
\chi_{3^2}(21^4) &=\chi_{31}(1^4)-\chi_{2^2}(1^4)=3-2=1. 
\end{align*}
So the one-to-one correspondence is not nicely explained by Murnaghan's formula. 

\subsection{Main result}
Here we are going to state our first main result, which explains the re-appearance phenomenon above. 

Let $m>1$ be an integer and $k$ be an integer. 
We define a $\CC$-linear operator $\phi_{k}^{(m)}: F \mor F$ by
\begin{equation}
\phi_{k}^{(m)} ([x_1, x_2, \ldots]) = \left( \sum_{a_0, \ldots, a_{m-1}}  [a_{m-1}, a_{m-2}, \ldots, a_0, x_1, x_2, \ldots]\right) - m  \label{defphi}
\end{equation}
for a maya diagram $[x_1, x_2, \ldots]$, where the sum is over all $m$-tuples of integers $(a_0, \ldots, a_{m-1})$ with $a_i \equiv i \pmod m$ and $a_0+\cdots+a_{m-1}=(0+\cdots+(m-1))-km$. 
Here the $-m$ on the right-hand side is defined as a $\CC$-linear operator $-m : \tilde{F} \mor \tilde{F}$ by $[x_1, x_2, \ldots]-m=[x_1-m, x_2-m, \ldots]$ for a maya diagram $[x_1, x_2, \ldots]$, 
which was applied to the sum in order that the right-hand side to lie in $F$. 

Our first result can now be stated as follows: 
\begin{thm}
	Let $\lmb$ be a maya diagram of size $n$ and let $\mu$ be a partition of $n-m$
	which does not have any multiple of $m$ as its part. 
	Then we have 
	\begin{equation*}
		\chi_{\lmb}(\mu \cup (m)) = \chi_{-\phi(\lmb)}(\mu). 
	\end{equation*}
	Here $\mu \cup (m)$ means the partition obtained by appending a part $m$ to $\mu$, $\phi=\phi_1^{(m)}$, 
	and the notation $\chi_{\kappa}(\nu)$ has been extended for $\kappa \in F$ linearly in $\kappa$. 
	\label{mainthm}
\end{thm}

\begin{eg}
Let $\lmb=(4,2,2) \leftrightarrow [3,0,-1,-4,-5,-6,\ldots]$ and $m=2$. 
Then
\begin{align*}
\phi(\lmb) &= \left(\sum_{\substack{a_0 \equiv 0, a_1 \equiv 1 \pmod 2 \\ a_0+a_1=-1} } [a_1, a_0,3,0,-1,-4,-5,-6,\ldots]\right)-2 \\
 &= \bigl([1,-2,3,0,-1,-4,-5,-6,\ldots]+[-3,2,3,0,-1,-4,-5,-6,\ldots]\bigr)-2 \\
 &= \bigl([3,1,0,-1,-2,-4,-5,-6,\ldots]-[3,2,0,-1,-3,-4,-5,-6,\ldots]\bigr)-2 \\
 &= [1,-1,-2,-3,-4,-6,-7,-8,\ldots]-[1,0,-2,-3,-5,-6,-7,-8,\ldots] \\
&\leftrightarrow (2,1,1,1,1)-(2,2,1,1). 
\end{align*}
Thus we have $\chi_{42^2}(\mu \cup (2))=-\chi_{21^4}(\mu)+\chi_{2^21^2}(\mu)$ for
 a partition $\mu$ with all parts odd. 
 Note that this is different from the one given by Murnaghan's formula, i.e. $\chi_{42^2}(\mu \cup (2)) = \chi_{2^3}(\mu) - \chi_{41^2}(\mu)+\chi_{42}(\mu)$ (although the latter is valid for any partition $\mu$).
\end{eg}

It can be easily seen that Theorem \ref{mainthm}, with $m=2$ and $\mu=(1, \ldots, 1)$, really explains the phenomenon. 
In fact, it can be shown that if $m=2$ then the sum in $\phi(\lmb)$ has at most $r-1$ nonzero terms for $|\lmb| < 2r^2$. 
So with $r=2$, it can be seen that $\chi_{\lmb}(21^*)$ can be expressed 
in the form $\pm \chi_{\bar{\lmb}}(1^*)$ (or zero) for a single diagram $\bar{\lmb}$ if $|\lmb| \leq 7$. 

In order to prove the theorem, we introduce some objects in order to simplify calculations in the proof. 
For $r \in \ZZ$, define a $\CC$-linear operator $\tilde{b}_r : \tilde{F} \mor \tilde{F}$ by 
$\tilde{b}_r[x_1, x_2, \ldots] =  [r, x_1, x_2, \ldots]$. 
Then 
\begin{equation}
\phi_{k}^{(m)}(\lmb) = \left(\sum \tilde{b}_{a_{m-1}} \cdots \tilde{b}_{a_0} \lmb \right) - m \label{phi2b}
\end{equation}
for any maya diagram $\lmb$, where the sum runs over the same $m$-tuples $(a_0, \ldots, a_{m-1})$ as in the definition (\ref{defphi}) of $\phi_{k}^{(m)}$. 
Let $\tilde{\phi}_{k}^{(m)} = \sum \tilde{b}_{a_{m-1}} \cdots \tilde{b}_{a_0}$. 
Then $\phi_{k}^{(m)}(\lmb) = \tilde{\phi}_{k}^{(m)}(\lmb) - m$. 

It is easy to check that the operators $\tilde{b}_r$ satisfy $\tilde{b}_r\tilde{b}_s=-\tilde{b}_s\tilde{b}_r$ and $[A_l, \tilde{b}_r]=\tilde{b}_{r-l}$ for integers $r, s$ and a positive integer $l$.

\begin{prf}
We prove the following two things: 
(i) $\phi$ commutes with $A_l$ ($m \nmid l$), 
and (ii) $-\phi(\lmb) = A_m \lmb$ for maya diagrams $\lmb$ of size $m$. 
In fact, with (i) and (ii) we can show the theorem as
\[
\chi_{\lmb}(\mu \cup (m)) \emp = A_mA_{\mu}\lmb=-\phi(A_{\mu}{\lmb})=A_{\mu}(-\phi({\lmb}))=\chi_{-\phi(\lmb)}(\mu) \emp. 
\]

We have for every $l$, 
\begin{align*}
[A_l, \tilde{\phi}_{k}^{(m)}]&=\sum_{j, (a_0, \ldots, a_{m-1})} \tilde{b}_{a_{m-1}} \cdots [A_l, \tilde{b}_{a_j}] \cdots \tilde{b}_{a_0} \\
&=\sum_{j, (a_0, \ldots, a_{m-1})} \tilde{b}_{a_{m-1}} \cdots \tilde{b}_{a_j-l} \cdots \tilde{b}_{a_0}. 
\end{align*}
Let us now assume that $l$ is not divisible by $m$. 
Consider the summand $\tilde{b}_{a_{m-1}} \cdots \tilde{b}_{a_j-l} \cdots \tilde{b}_{a_0}$ indexed by $(j, (a_0, \ldots, a_{m-1}))$. 
Take $j' \in \{0, \ldots, m-1\}$ such that $j' \equiv j-l \pmod m$. 
If $a_j-l = a_{j'}$, then the summand is zero. 
Otherwise, it cancels out with another summand indexed by $(j, (a'_0, \ldots, a'_{m-1}))$, 
where $a'_j=a_{j'}+l$, $a'_{j'}=a_j-l$ and $a'_{j''}=a_{j''}$ for $j'' \neq j, j'$ 
(note that the above correspondence $(j, (a_0, \ldots, a_{m-1})) \mapsto (j, (a'_0, \ldots, a'_{m-1}))$ is involutive). 
Thus we have $[A_l, \tilde{\phi}_{k}^{(m)}]=0$ and this shows (i) (it is clear that $A_l$ commutes with shifting). 

Let us now show (ii). It is well known that, for a partition $\lmb$ of size $m$, 
$A_m \lmb = \chi_\lmb((m))\emp = (-1)^{r-1}\emp$ if $\lmb$ is a hook with $r$ rows and $0$ otherwise. 
In our maya-diagram setting, this may be formulated as
\begin{equation}
A_m \lmb = \left\{
\begin{array}{ll}
(-1)^{r-1}\emp & (\lmb = [-r+m, -1, -2, \ldots, \widehat{-r}, \ldots, -m, -m-1, \ldots]) \\
0 & (\text{otherwise})
\end{array}\right. \label{mn_hook}
\end{equation}
where $\widehat{-r}$ means that $-r$ is to be removed from the sequence. 
Thus it suffices to show that $-\phi(\lmb)$ coincides with the right-hand side of (\ref{mn_hook}). 

Assume now $\phi(\lmb) \neq 0$. 
Let $\lmb = [x_1, x_2, \ldots]$ and let $y_1<y_2<\cdots$ be the integers not appearing in $x_1, x_2, \ldots$. 
By Lemma \ref{mayalem}, we have $y_i = i-1$ for $i > m$ and $\sum (i-1-y_i) = m$. 
On the other hand, 
from $\phi(\lmb)\neq 0$ we have $[a_{m-1}, a_{m-2}, \ldots, a_0, x_1, x_2, \ldots] \neq 0$
for some $(a_0, \ldots, a_{m-1})$ with $a_i \equiv i \pmod m$ and $a_0+\cdots+a_{m-1} = (0+\cdots+(m-1))-m$. 
$[a_{m-1}, a_{m-2}, \ldots, a_0, x_1, x_2, \ldots] \neq 0$ implies that all terms $a_i$ appear in the sequence $(y_i)$, but
\begin{align*}
a_0+\cdots+a_{m-1} &= (0+\cdots+m-1)-m \\
&= (0+\cdots+m-1)-\sum_{i=1}^\infty (i-1-y_i) \\
&= (0+\cdots+m-1)-\sum_{i=1}^m (i-1-y_i) \\
&= y_1+\ldots+y_m
\end{align*}
so we have $\{ y_1, \ldots, y_m \} = \{ a_0, \ldots, a_{m-1} \} \equiv \{ 0, \ldots, m-1 \} \pmod m$. 
Combining $y_1+\ldots+y_m=(0+\cdots+m-1)-m$, $\{ y_1, \ldots, y_m \} \equiv \{ 0, \ldots, m-1 \} \pmod m$ and $y_1 < \ldots < y_m < y_{m+1} = m$, 
we can conclude $(y_1, \ldots, y_m) = (s-m, 0, 1, \ldots, \widehat{s}, \ldots, m-1)$ for some $0 \leq s \leq m-1$. 
This means $\lmb = [-r+m, -1, -2, \ldots, \widehat{-r}, \ldots, -m, -m-1, \ldots]$ with $r=m-s$, 
and $\phi(\lmb)=(-1)^r\emp$ can be easily checked in these cases (recall that we have only one nonzero term in $\phi(\lmb)$). 
\hfill $\Box$
\end{prf}

Consider the case $m \mid l$ in (i) above. 
The same calculations as in the proof of (i) above yields $[A_{rm}, \phi_k^{(m)}]=m\phi_{k+r}^{(m)}$. 
Using this commutation relation, one can yield more general expansion of the character values
$\chi_\lmb(\mu \cup \nu)$ into the values $\chi_\kap(\mu)$, where $\mu$ and $\nu$ are partitions, and no parts of $\mu$ and all parts of $\nu$ are divisible by $m$. 
For example, for $\nu=(m^k)$ one has the following: 
\begin{thm}
For a partition $\mu$ with no parts divisible by $m$ one has
\begin{equation}
\chi_\lmb(\mu \cup m^k) = -\sum_{i=0}^{k-1} (-m)^i \binom{k-1}{i} \chi_{\phi_{i+1}^{(m)}(\lmb)}(\mu \cup m^{k-1-i}). \label{multim}
\end{equation}
\end{thm}
\begin{prf}
$[A_{m}, \phi_k^{(m)}]=m\phi_{k+1}^{(m)}$ yields
$\phi_1^{(m)}A_{m}^k=A_{m}^k\phi_1^{(m)}+\sum_{i=0}^{k-1} (-m)^i \binom{k-1}{i} \phi_{i+1}^{(m)}A_m^{k-1-i}$, 
and the conclusion easily follows from this by the same argument as in the proof of Theorem \ref{mainthm} above. 
\hfill $\Box$
\end{prf}

The case $\mu=\emp$, $\nu=(1^k)$ and $m=1$ is particularly interesting. 
In this case, the formula gives an expansion for the degree of $\chi_\lmb$ 
in a fashion similar to the ordinary determinant formula
\[
\chi_\lmb(1^k) = \sum_{w \in S_l} \mathrm{sgn}(w)g(w(\lmb+\rho)-\rho)
\text{ where }g(\alpha)=(\sum_i \alpha_i)! / \prod_i (\alpha_i!)
\]
for the degree, though different from it: 
\begin{thm}
For a partition $\lmb$ with size $k$ and $l$ rows, 
one has
\[
\chi_\lmb(1^k) = \sum_{w \in S_l} \mathrm{sgn}(w)f(w(\lmb+\rho)-\rho)
\]
where $\rho=(l-1, \ldots, 0)$. 
Here for $\alpha=(\alpha_1,\cdots,\alpha_l) \in \ZZ^l$, 
we set $f(\alpha)=\prod_{i=1}^{l} \binom{\alpha_i+\cdots+\alpha_l-1}{\alpha_i-1}$ if $\alpha_i > 0$ for all $i$
 and $f(\alpha)=0$ if $\alpha_i \leq 0$ for some $i$. 
\end{thm}

\begin{prf}
We show $\chi_{\lmb'}(1^k) = \sum_{w \in S_l} \mathrm{sgn}(w)f(w(\lmb+\rho)-\rho)$, where $\lmb'$ is the conjugate partition of $\lmb$. 
Let $[x_1, x_2, \ldots]$ and $[y_1, y_2, \ldots]$ be the maya diagrams corresponding to $\lmb$ and $\lmb'$ respectively. 

Let $\phi_i=\phi_i^{(1)}$. 
Successive application of (\ref{multim}) gives
\begin{align*}
\chi_{\lmb'}(1^k) \emp &= \sum (-1)^{a_1+\ldots+a_r} \left( \prod_{i=1}^r \binom{k-a_1-\cdots-a_{i-1}-1}{a_i-1}\right) \phi_{a_r} \cdots \phi_{a_1}(\lmb') \\
&= (-1)^k \sum \left( \prod_{i=1}^r \binom{a_i+\cdots+a_r-1}{a_i-1}\right) \phi_{a_r} \cdots \phi_{a_1}(\lmb')
\end{align*}
where the sum is over all $r \geq 0$ and $a_1, \ldots, a_r \in \ZZ_{>0}$ such that $a_1+\ldots+a_r=k$. 

In general, 
if $i>0$ and $\kap$ is a partition with $c$ columns, 
it can be shown that $\phi_i(\kap)$ is either $\pm 1$ times a partition with exactly $c-1$ columns or zero. 
Thus if $\phi_{a_r} \cdots \phi_{a_1}(\lmb')$ is a nonzero multiple of $\emp$, one must have $r=l$. 

Since $\phi_{a_l} \cdots \phi_{a_1}(\lmb') = [-a_l-1, \ldots, -a_1-l, y_1-l, y_2-l, \ldots]$, 
in order for this to be nonzero the numbers $-a_i+i-1$ must lie in the complement of $\{ y_1, y_2, \ldots \}$, 
which is, by Lemma \ref{mayalem}, $\{ -1-x_1, -1-x_2, \ldots \}$. 
In other words, all $a_i-i$ must be in $\{ x_1, x_2, \ldots \}$. 
Moreover, we have $a_i-i \geq -i > -l-1 = x_{l+1}$ for $1 \leq i \leq l$ so $a_i-i \in \{x_1, \ldots, x_l\}$. 
And since $-a_l-1, \ldots, -a_1-l$ must be all distinct, all $a_i-i$ must be distinct. 
This shows $a_1-1, \ldots, a_l-l$ must be a permutation of $x_1, \ldots, x_l$. 
That is, $a+\rho=w(\lmb+\rho)$ for a permutation $w \in S_l$.
Then $
\phi_{a_l} \cdots \phi_{a_1}(\lmb')=\mathrm{sgn}(w)[-1-x_l-l, \ldots, -1-x_1-l, y_1-l, y_2-l, \ldots]$. 
It can be easily shown that $[-1-x_l-l, \ldots, -1-x_1-l, y_1-l, y_2-l, \ldots]=(-1)^k\emp$ and this completes the proof. 
\hfill $\Box$
\end{prf}


\section{Relation with Bernstein operator}
Let $\Lmb$ denote the ring of symmetric functions in the variables $X_1, X_2, \ldots$: 
$\Lmb=\ZZ[e_1,e_2,\ldots]=\ZZ[h_1,h_2,\ldots]$ where $e_r=\sum_{i_1<\ldots<i_r} X_{i_1} \cdots X_{i_r}$ and $h_r=\sum_{i_1 \leq \ldots \leq i_r} X_{i_1} \cdots X_{i_r}$. 
We define $h_0=e_0=1$ and $h_r=e_r=0$ for $r<0$. 
Let $H(u):=\sum_{n \geq 0} h_nu^n = \prod_{i \geq 1} \frac{1}{1-X_iu}$ and $E(u):=\sum_{n \geq 0} e_nu^n = \prod_{i \geq 1} (1+X_iu)$ be their generating functions. 
Let $p_r=\sum_{i \geq 1} X_i^r$ ($r \geq 1$). 
We denote the Schur function by $s_\lmb$ (\cite[\S I.3]{Mac}). 
In this paper, the index of a Schur function may be either a partition, written in parentheses, or a maya diagram, written in brackets.

We consider on $\Lmb$, as usual, the Hall inner product $\left< \cdot , \cdot \right> : \Lmb \times \Lmb \mor \ZZ$ by 
$\left< s_\lmb, s_\mu \right> = \delta_{\lmb\mu}$ for partitions $\lmb, \mu$. 
Let $f^\perp$ denote the adjoint operator of the multiplication by $f$ with respect to this inner product, 
say: $\left< f^\perp(g), h \right> = \left< g, fh\right>$. 

Using the Schur functions and the Hall inner product, 
Theorem \ref{mnformula} can be stated as $p_{\mu_1}^\perp \cdots p_{\mu_l}^\perp s_\lmb = \chi_{\lmb}(\mu)$ 
(\cite[\S I.7]{Mac}), which is in fact much closer to the original Frobenius formula. 

Bernstein's creation operator is defined as follows: 
\begin{defn}
For $n \in \ZZ$, 
\begin{equation}
B_n = \sum_{i \geq 0} (-1)^i h_{n+i} e_i ^ \perp \label{berndef}
\end{equation}
where $h_{n+i}$ denotes the multiplication operator by $h_{n+i}$. 
Note that, even though the expression above is an infinite sum, only finitely many terms give nonzero image when applied to each $f \in \Lmb$. 
\end{defn}

Bernstein operators have the following property (\cite[\S I.5, Example 29]{Mac}): 
\begin{prop}
For a partition $\lmb = (\lmb_1, \ldots, \lmb_l)$, 
\[
B_{\lmb_1} \cdots B_{\lmb_l}(1) = s_{\lmb}. 
\]
\end{prop}
This is still valid for a general integer sequence $\lmb$
if we interpret Schur functions indexed by integer sequences which is nondecreasing or having nonpositive parts 
as $s_{(\cdots,i,j,\cdots)}=-s_{(\cdots,j-1,i+1,\cdots)}$. 
In our maya-diagram setting, this implies
\begin{equation}
B_n(s_{[x_1, x_2, \ldots]})=s_{[n, x_1, x_2, \ldots]-1}, 
\label{bmaya}
\end{equation}
or, for $\lmb$ a maya diagram, 
\[
B_n(s_{\lmb})=s_{b_n\lmb}, 
\]
where $b_n\lmb=\tilde{b}_n\lmb-1$ with $\tilde{b}_n$ defined just before the proof of theorem \ref{mainthm}. 
From this we can see $B_{a_{m-1}} \cdots B_{a_0}s_{[x_1, x_2, \ldots]}=s_{[a_{m-1}-1, \ldots, a_0-m, x_1-m, x_2-m, \ldots]}=s_{[a_{m-1}+m-1, \ldots, a_0, x_1, x_2, \ldots]-m}$, 
and so if we identify $F$ with $\Lmb_\CC=\Lmb \otimes \CC$ by identifying each maya diagram $\lmb$ with the Schur function $s_\lmb$ we get
\begin{equation}
\phi_{k}^{(m)} = \sum B_{a_{m-1}} \cdots B_{a_0}
\label{phiinb}
\end{equation}
where sum runs over all m-tuples $(a_0, \ldots, a_{m-1})$ with $a_i \equiv 0 \pmod m$ and $a_0+\cdots+a_{m-1}=-km$. 

We have the following Bernstein-operator like expression for $\phi_{k}^{(m)}$: 
\begin{thm}
For $m \geq 1$ and $n \in \ZZ$ we have
\begin{equation}
\phi_{-n}^{(m)} = \sum_{i \geq 0} (-1)^i (h_{n+i} \circ p_m) (e_i \circ p_m) ^ \perp
\label{bernpm}
\end{equation}
where $- \circ p_m$ denotes the plethysm with the $m$-th power sum: $(f \circ p_m)(\{X_i\})=f(\{X_i^m\})$. 
\label{bernpm_thm}
\end{thm}

\begin{rem}
In fact, the properties (i) and (ii) of $\phi$ in the proof of Theorem \ref{mainthm} can be also seen from the above expression for $\phi_{n}^{(m)}$. 
\end{rem}

\begin{prf}
Let $E^\perp(u)=\sum_{n \geq 0} e_n^\perp u^n$ and 
$B(u)=\sum_{n \in \ZZ} B_nu^n = H(u)E^\perp(-u^{-1})$. 
It is known that $B(u)$ satisfies the following (\cite[\S I.5, Example 29]{Mac}): 
\begin{equation}
B(u_1)\cdots B(u_r)=\prod_{p<q}(1-u_p^{-1}u_q) \cdot H(u_1)\cdots H(u_r) E^\perp(-u_1^{-1}) \cdots E^\perp(-u_r^{-1}).
\label{commB}
\end{equation}
Let $(j_0, \ldots, j_{m-1}) \in (\ZZ/m\ZZ)^m$ and $\omega$ be a primitive $m$-th root of unity. 
Letting $(u_1, \ldots, u_r)=(\omega^{j_0}u, \ldots, \omega^{j_{m-1}}u)$ in (\ref{commB}), we have
\begin{equation}
\begin{split}
&B(\omega^{j_0}u) \cdots B(\omega^{j_{m-1}}u) \\
& \; = \prod_{p<q}(1-\omega^{j_q-j_p}) \\
& \qquad \cdot H(\omega^{j_0}u)\cdots H(\omega^{j_{m-1}}u) E^\perp(-\omega^{-j_0}u^{-1}) \cdots E^\perp(-\omega^{-j_{m-1}}u^{-1}). 
\end{split}
\label{hoge}
\end{equation}
Here, the well-definedness of the left-hand side of (\ref{hoge}) follows from the following claim, which can be easily seen by using (\ref{bmaya}): 
\begin{quote}
\noindent\textbf{Claim. }
For any $s \geq 0$ and $f \in \Lmb$, if $n \in \ZZ$ is sufficiently small 
then for any $a_1, \ldots, a_s \in \ZZ$ we have $B_nB_{a_1} \ldots B_{a_s}(f)=0$. 
\hfill $\Box$
\end{quote}
The right-hand side of (\ref{hoge}) is zero if some $j_p, j_q$ ($p \neq q$) are equal, 
and otherwise it is a constant multiple of $H(u) \cdots H(\omega^{m-1}u) E^\perp(-u^{-1}) \cdots E^\perp(-\omega^{m-1}u^{-1})$. 
Thus, summing (\ref{hoge}) over all $(j_0, \ldots, j_{m-1})$ and dividing by $m^m$, we have
\[
\left( \sum_{n \in \ZZ} B_{nm}u^{nm}\right)^m
= C_m \cdot H(u) \cdots H(\omega^{m-1}u) E^\perp(-u^{-1}) \cdots E^\perp(-\omega^{m-1}u^{-1}), 
\] 
where $C_m$ is a constant depending only on $m$. 
Notice that the left-hand side is equal to $\sum_n \phi_{-n}^{(m)} u^{nm}$. 

We have
$\displaystyle H(u)\cdots H(\omega^{m-1}u) = \prod_{\substack{i \geq 1 \\ j \in \ZZ/m\ZZ}} \frac{1}{1-\omega^jX_iu} = \prod_{i \geq 1} \frac{1}{1-X_i^mu^m} = \sum_{n \geq 0} (h_n \circ p_m) u^{nm}$. 
We also have 
$\displaystyle E(-u^{-1})\cdots E(-\omega^{m-1}u^{-1}) = \prod_{\substack{i \geq 1 \\ j \in \ZZ/m\ZZ}} (1-\omega^jX_iu^{-1}) = \prod_{i \geq 1} (1-X_i^mu^{-m}) = \sum_{n \geq 0} (-1)^n(e_n \circ p_m) u^{-nm}$ 
and thus 
$E^\perp(-u^{-1})\cdots E^\perp(-\omega^{m-1}u^{-1}) = \sum_{n \geq 0} (-1)^n(e_n \circ p_m)^\perp u^{-nm}$. 
Thus we have $\phi_{-n}^{(m)} = C_m \sum_{i \geq 0} (-1)^i (h_{n+i} \circ p_m) (e_i \circ p_m) ^ \perp$. 
Setting $n=0$ and comparing the actions of both side on $1$ shows $C_m=1$. 
\hfill $\Box$
\end{prf}

\section{An analog for projective characters}
A \textit{projective representation} $\pi=(V, \pi)$ of a group $G$ is a group homomorphism $\pi$ from $G$ to $PGL(V)$, 
the projective general linear group of a vector space $V$. 
Two projective representations $(V, \pi)$ and $(W, \rho)$ are said to be \textit{projectively equivalent} 
if there exists a vector space isomorphism $V \mor W$ such that 
the induced group isomorphism $f: PGL(V) \mor PGL(W)$ satisfies $f\pi(g)=\rho(g)f$ for all $g \in G$.
Let us call a projective representation $\pi$ \textit{nontrivial} 
if $\pi$ is not projectively equivalent to any representations obtained from a linear representation of $G$ 
by composing with $GL(V) \surj PGL(V)$. 

Let $\tilde{S}_n$ be the group generated by the generators $s_1, \ldots, s_{n-1}, z$ bound by the relations: 
\begin{itemize}
\item $z$ is a central element with $z^2=1$, 
\item $s_i^2=z$, $s_is_j=zs_js_i$ ($|i-j| \geq 2$), $(s_is_{i+1})^3=z$. 
\end{itemize}

Clearly one has a surjective group homomorphism $\theta: \tilde{S}_n \surj S_n$ which sends $z$ to $1$ and $s_i$ to $(i\;i+1)$, $i=1, \ldots, n-1$. 
This homomorphism has kernel $\{ 1,z \}$. 
Since $z$ is in the center of $\tilde{S}_n$ and has order 2, its action on an irreducible representation is either by $1$ or by $-1$. 
Call an irreducible representation of $\tilde{S}_n$ \textit{negative} if $z$ acts by $-1$. 
Call two representations of $\tilde{S}_n$ being \textit{associate} of each other 
if one can be obtained from the other by tensoring with $sgn \circ \theta$, where $sgn$ is the sign representation of $S_n$. 
The following relationship between the projective representations of $S_n$ and the linear representations of $\tilde{S}_n$ is known: 
\begin{prop}[{\cite[Chap. 2]{HH}}]
If $n \geq 4$, Projective isomorphism classes of nontrivial irreducible projective representations of $S_n$
 is in one-to-one correspondence with the associate classes of negative irreducible representations of $\tilde{S}_n$. 
\end{prop}

It is known that the isomorphism classes of negative representations of $\tilde{S}_n$ are indexed by the strict partitions of $n$ (a partition $(\lmb_1,\ldots,\lmb_l)$ is called \textit{strict} if $\lmb_1>\cdots>\lmb_l$). 
Let $\psi_\lmb$ denote the irreducible character of $\tilde{S}_n$ indexed by $\lmb$. 

If $C$ is a conjugacy class of $S_n$, 
$\theta^{-1}(C)$ is either a single conjugacy class or the union of two conjugacy classes. 
In the former case, $g$ and $zg$ are conjugate for $g \in \theta^{-1}(C)$ 
and thus negative irreducible characters vanish there. 
In the latter case we say that $C$ \textit{splits}. 
Conjugacy class $C_\mu$ of $S_n$ with cycle type $\mu$ splits iff: 
(i) all parts of $\mu$ are odd (in which case we call $\mu$ \textit{all-odd}), 
or (ii) $\mu$ is strict and $C_\mu$ consists of odd permutations (\cite[Theorem 3.8]{HH}). 
In fact it is easy to describe the character values explicitly in the case (ii), 
so we are interested in the case (i). 
Let $\psi_\lmb(\mu)$ denote the value of $\psi_\lmb$ evaluated at an element $g_\mu$, 
which is chosen from $\theta^{-1}(C_\mu)$ 
so that the character of the ``basic representation'' (\cite[Chap. 6]{HH}) of $\tilde{S}_n$ 
takes a positive value at $g_\mu$. 
We also let $\tilde{\psi}_\lmb(\mu) = 2^{\lceil \frac{\ell(\lmb)-\ell(\mu)}{2} \rceil}\psi_\lmb(\mu)$. 

Let $q_n=\sum h_{n-i}e_i \in \Lmb$ and let $\Gamma$ be the subring of $\Lmb$ generated by $q_1, q_2, q_3, \ldots$. 
It is known that $\Gamma\otimes \QQ = \QQ[p_1, p_3, p_5, p_7, \ldots]$. 
We have a basis $\{Q_\lmb\}_{\text{$\lmb$:strict partition}}$ of $\Gamma$
consisting of so-called \textit{Schur $Q$-functions} (\cite[Chap. 7]{HH}, \cite[\S III.8]{Mac}). 
Just as Schur functions carry information about irreducible characters of linear representations of symmetric groups, 
Schur $Q$-functions carry information about projective characters $\psi_\lmb$: 
in fact, if we define an inner product $\left< , \right>$ on $\Gamma$ 
by $\left< Q_\lmb, Q_\mu \right> = \delta_{\lmb\mu}2^{\ell(\lmb)}$ for all strict partitions $\lmb$ and $\mu$ 
and denote the adjoint of the multiplication by $f\in \Gamma$ as $f^\perp$, 
then for strict $\lmb$ and all-odd $\mu$ we have 
$\tilde{\psi}_\lmb(\mu)=p_{\mu_1}^\perp \cdots p_{\mu_l}^\perp Q_\lmb$ (\cite[Chap. 8]{HH}). 

As in the Schur-funtion case, Schur $Q$-functions also have creation operators: 
\begin{prop}[{\cite[Theorem 7.21]{HH}}]
If we define $\displaystyle \QB_n = \sum_{i \geq 0} (-1)^i q_{n+i}q_{i}^\perp$ for $n \in \ZZ$, 
then for a strict partition $\lmb$ we have $Q_\lmb=\QB_{\lmb_1} \cdots \QB_{\lmb_l}(1)$. 
\end{prop}
Let $Q_\alpha = \QB_{\alpha_1}\cdots\QB_{\alpha_r}(1)$ for all $\alpha=(\alpha_1, \ldots, \alpha_r) \in \ZZ^r$. 
The operators $\QB_r$ satisfy $\QB_r\QB_s+\QB_s\QB_r=2(-1)^r\delta_{r,-s}$ (\cite[Theorem 9.1]{HH}), 
and from this one has ``reordering rules'' for writing $Q_\alpha$ as a linear combination of the functions $Q_\lmb$ with strict partitions $\lmb$: 
\begin{itemize}
\item If, for some $i \geq 1$, the subsequence of $\alpha$ consisting of all occurrences of $\pm i$ is \textit{not} of the form
$i, -i, i, \ldots, -i, i$ or $-i, i, \ldots, -i, i$, then $Q_\alpha=0$. 
\item Otherwise, there exists a permutation of the sequence $\alpha$ 
which has the form $\lmb, -a_1, a_1, \ldots, -a_r, a_r, 0, \ldots, 0$
for a strict partition $\lmb$ and positive integers $a_1, \ldots, a_r$ (not necessarily distinct). 
In this case $Q_\alpha=(-1)^{a_1+\ldots+a_r}2^r\epsilon Q_\lmb$, 
where $\epsilon$ is the sign of any permutation which permutes $\alpha$ into the form above
while, for each $i \geq 0$, keeping the order of the terms $\pm i\,$. 
\end{itemize}

Since the definition of $\QB_n$ is similar to the definition (\ref{berndef}) of the Bernstein operator, 
we can consider a modification of $\QB_n$ analogous to (\ref{bernpm}): 
let $\Phi_{n}^{(m)}=\sum_{i \geq 0} (-1)^i (q_{n+i} \circ p_m) (q_{i} \circ p_m)^\perp : \Gamma \mor \Gamma$ for $m\geq 1$ odd and $n \in \ZZ$. 
Then we have the following formula for $\Phi_{n}^{(m)}$ analogous to (\ref{phiinb}): 
\begin{thm}
For any $m \geq 1$ odd, we have
\begin{align*}
&\sum_n \Phi_{n}^{(m)} u^{nm}\\
&= 
\left( 1+\sum_{\substack{i,j \in \ZZ \\ i \equiv -2 \\ j \equiv 2}} \QB_{i,j}u^{i+j} \right)
\left( 1+\sum_{\substack{i,j \in \ZZ \\ i \equiv -4 \\ j \equiv 4}} \QB_{i,j}u^{i+j} \right) \cdots
\left( 1+\sum_{\substack{i,j \in \ZZ \\ i \equiv -m+1 \\ j \equiv m-1}} \QB_{i,j}u^{i+j} \right)
\left( \sum_{\substack{k \in \ZZ \\ k \equiv 0}} \QB_ku^k \right)
\end{align*}
where the congruences are modulo $m$
and $\QB_{i,j} = \left\{ \begin{array}{ll} \QB_i\QB_j & (i>j) \\ -\QB_j\QB_i & (i<j) \\ 0 & (i=j) \end{array} \right.$. 
\label{projthm}
\end{thm}
We note that in the product above the coefficient of each $u^d$ is well-defined 
by the reordering rule above. 

Since $\Phi_{-1}^{(m)}$ commutes with $p_l^\perp$ ($m \nmid l$) and coincides with a constant multiple of $p_m^\perp$ (in fact $-2 p_m^\perp$) on degree $m$ part of $\Gamma$ 
by the same reason as the remark after Theorem \ref{bernpm_thm}, 
we have a corresponding relation for $\tilde{\psi}_\lmb(\mu)$ as in Theorem \ref{mainthm}:
\begin{cor}
Let $\mathcal{F}$ be the vector space generated by all integer sequences of finite length
bound by the same relations as the reordering rules for the $Q$-functions
(so a sequence $\alpha$ is equal to zero in $\mathcal{F}$ if there exists $i \geq 1$ such that the subsequence of $\alpha$ consisting of all occurences
of $\pm i$ is not of the form $i, -i, i, \ldots, -i, i$ or $-i, i, \ldots, -i, i$, 
and equals to some nonzero constant multiple of a strict partition otherwise). 
Then for an odd integer $m \geq 1$, a strict partition $\lmb \in \mathcal{F}$ and an all-odd partition $\mu$ with no parts divisible by $m$ such that $|\lmb|=|\mu|+m$, 
one has $\tilde{\psi}_\lmb(\mu \cup (m))=\tilde{\psi}_{-\frac{1}{2}\Phi(\lmb)}(\mu)$. 
Here $\Phi: \mathcal{F} \mor \mathcal{F}$ is a $\CC$-linear map 
defined by 
$\Phi(\lmb)= \sum_\alpha \epsilon_\alpha\cdot(\alpha, \lmb)$, 
where the sum is over all integer sequences $\alpha=(\alpha_1, \ldots, \alpha_{2s+1})$ of odd length
such that
$\alpha_{2i-1}>\alpha_{2i}$, $\{\alpha_{2i-1}, \alpha_{2i}\} \equiv \{2c_i, -2c_i\} \pmod{m}$
for some integers $\frac{m-1}{2}\geq c_1>\cdots>c_s>0$ and $\alpha_{2s+1} \equiv 0 \pmod{m}$
and $\sum_i \alpha_i=-m$. 
For such a sequence, 
$\epsilon_\alpha=(-1)^{n(\alpha)}$ where $n(\alpha)$ is the number of $i$ such that $\alpha_{2i} \equiv -2c_i \pmod{m}$
and $(\alpha,\lmb)$ is the sequence obtained by concatenating $\alpha$ and $\lmb$. 

\hfill $\Box$
\end{cor}

\begin{eg}
\begin{align*}
\Phi_{-1}^{(3)}Q_{4,3,2}&=Q_{\underline{1,-4,0},4,3,2}+Q_{\underline{4,-4,-3},4,3,2}-Q_{\underline{2,-2,-3},4,3,2}+Q_{\underline{-3},4,3,2} \\
&=-2Q_{3,2,1}+4Q_{4,2}-4Q_{4,2}+2Q_{4,2} \\
&=-2Q_{3,2,1}+2Q_{4,2}, 
\end{align*}
so we have $\tilde{\psi}_{4,3,2}(\mu \cup (3)) = \tilde{\psi}_{3,2,1}(\mu)-\tilde{\psi}_{4,2}(\mu)$ for $\mu$ with no parts divisible by $3$, say $\mu=(1^6)$ and $(5,1)$. 
\end{eg}

\noindent{\textbf{Proof of Theorem \ref{projthm}. }}
Let $Q(u)=\sum_{n \geq 0} q_nu^n$, $Q^\perp(u)=\sum_{n \geq 0} q_n^\perp u^n$ and $\QB(u)=\sum_{n \in \ZZ} \QB_nu^n=Q(u)Q^\perp(-u^{-1})$. 
Let $\omega$ be a primitive $m$-th root of unity
and define $\displaystyle \QB(a,b;u) = \frac{1-\omega^{b-a}}{1+\omega^{b-a}}+\sum_{i,j \in \ZZ} \omega^{ai+bj}\QB_{i,j}u^{i+j}$ for $a,b \in \ZZ/m\ZZ$. 
\begin{tmplem}
We have
\begin{equation}
\QB(a,b;u) = \frac{1-\omega^{b-a}}{1+\omega^{b-a}} Q(\omega^a u)Q(\omega^b u)Q^\perp(-\omega^{-a}u^{-1})Q^\perp(-\omega^{-b}u^{-1}). 
\label{qbabu}
\end{equation}
\end{tmplem}
\begin{prf}
Let $\mathcal{C}_{r,s} = \sum_{i,j \geq 0} (-1)^{i+j}q_{r+i}q_{s+j}q_i^\perp q_j^\perp$ ($=\mathcal{C}_{s,r}$) for $r, s \in \ZZ$.  
It is easy to see that the right-hand side of (\ref{qbabu}) is equal to 
$\displaystyle \frac{1-\omega^{b-a}}{1+\omega^{b-a}}\sum_{r,s \in \ZZ} \omega^{ar+bs}\mathcal{C}_{r,s} u^{r+s}$. 
On the other hand, we have (\cite[Chap. 9, (7)]{HH}) $q_r^\perp q_s = q_s q_r^\perp+2\sum_{i \geq 1} q_{s-i}q_{r-i}^\perp$
and thus $\QB_r\QB_s=\mathcal{C}_{r,s}+2\sum_{i \geq 1} (-1)^i\mathcal{C}_{r-i,s-i}$
(in particular, we have $\mathcal{C}_{r,r}+2\sum_{i \geq 1} (-1)^i\mathcal{C}_{r-i,r-i} = \QB_r\QB_r = \delta_{r,0}$). 
Using these relations to rewrite both sides into linear combinations of $\mathcal{C}_{r,s}$ ($r>s$) and comparing the coefficients implies the Lemma. 
\hfill $\Box$
\end{prf}
Let $l \in \{2,4,\ldots,m-1 \}$. Multiplying (\ref{qbabu}) by $\omega^{l(a-b)}$, summing over $a,b \in \ZZ/m\ZZ$ and dividing by $m^2$ we have 
\begin{align*}
&\frac{1}{m^2}\sum_{a,b \in \ZZ/m\ZZ} \omega^{l(a-b)}\frac{1-\omega^{b-a}}{1+\omega^{b-a}} + \sum_{\substack{i,j \in \ZZ \\ i \equiv -l, j \equiv l \pmod m}} \QB_{i,j}u^{i+j}\\
&=\frac{1}{m^2}\sum_{a,b \in \ZZ/m\ZZ} \omega^{l(a-b)}\frac{1-\omega^{b-a}}{1+\omega^{b-a}} Q(\omega^a u)Q(\omega^b u)Q^\perp(-\omega^{-a}u^{-1})Q^\perp(-\omega^{-b}u^{-1}). 
\end{align*}
The first term on the left-hand side can be easily shown to be $1$. 
Since $q_r^\perp q_s = q_s q_r^\perp+2\sum_{i \geq 1} q_{s-i}q_{r-i}^\perp$ we have $Q^\perp(s)Q(t)=F(st)Q(t)Q^\perp(s)$
where $F(z)=1+2z+2z^2+2z^3+\cdots=\frac{1+z}{1-z}$, and thus
\begin{equation}
\begin{split}
&\left( 1+\sum_{\substack{i,j \in \ZZ \\ i \equiv -2 \\ j \equiv 2}} \QB_{i,j}u^{i+j} \right)
\left( 1+\sum_{\substack{i,j \in \ZZ \\ i \equiv -4 \\ j \equiv 4}} \QB_{i,j}u^{i+j} \right) \cdots
\left( 1+\sum_{\substack{i,j \in \ZZ \\ i \equiv -m+1 \\ j \equiv m-1}} \QB_{i,j}u^{i+j} \right)
\left( \sum_{\substack{k \in \ZZ \\ k \equiv 0}} \QB_ku^k \right) \\
&=\frac{1}{m^m}\sum_{\substack{a_1, b_1, \ldots, a_r, b_r, c \\ \in \ZZ/m\ZZ}} 
\left( \begin{array}{l}
\omega^{\sum_{i=1}^{r} 2i(a_i-b_i)} \\
\times \prod^{\longrightarrow}_{i=1,2,\ldots,r}
\left(
\begin{array}{l}
\frac{1-\omega^{b_i-a_i}}{1+\omega^{b_i-a_i}} Q(\omega^{a_i} u)Q(\omega^{b_i} u) \\
\qquad \cdot Q^\perp(-\omega^{-a_i}u^{-1})Q^\perp(-\omega^{-b_i}u^{-1})
\end{array}\right) \\
\times Q(\omega^c u)Q^\perp(-\omega^{-c} u^{-1})
\end{array} \right) \\
&=\frac{1}{m^m}\sum_{\substack{a_1, b_1, \ldots, a_r, b_r, c \\ \in \ZZ/m\ZZ}} 
\left( \begin{array}{l}
\omega^{\sum_{i=1}^{r} 2i(a_i-b_i)}
\prod_{1 \leq i < j \leq m} \frac{1-\omega^{c_j-c_i}}{1+\omega^{c_j-c_i}} \\
\times Q(\omega^{c_1}u) \cdots Q(\omega^{c_m}u) Q^\perp(-\omega^{-c_1}u^{-1}) \cdots Q^\perp(-\omega^{-c_m}u^{-1})
\end{array} \right)
\end{split}
\label{chounagai}
\end{equation}
where we set $r=\frac{m-1}{2}$ and $(c_1, \ldots, c_m)=(a_1, b_1, \ldots, a_r, b_r, c)$, and the congruences on the leftmost side are modulo $m$. 
We note that the well-definedness of the leftmost side follows from the following claim, which can be seen from the reordering rule above: 
\begin{quote}
\noindent\textbf{Claim. }
For any $s \geq 0$ and $f \in \Gamma$, if $n \in \ZZ$ is sufficiently small 
then for any $a_1, b_1, \ldots, a_s, b_s, d, i\in \ZZ$ 
we have $\QB_{i,n-i}\QB_{a_1,b_1} \ldots \QB_{a_s,b_s}\QB_d(f)=0$. 
\hfill $\Box$
\end{quote}

The summand of the rightmost side of (\ref{chounagai}) vanishes if some $c_i, c_j$ ($i \neq j$) are equal, 
and otherwise it is a constant multiple of $Q(u) \cdots Q(\omega^{m-1}u) \cdot Q^\perp(-u^{-1}) \cdots Q^\perp(-\omega^{m-1}u^{-1})$. 
Since $Q(u)=\prod_{i \geq 1} \frac{1+X_iu}{1-X_iu}$, the same calculation as in the proof of Theorem \ref{bernpm_thm} yields 
$Q(u) \cdots Q(\omega^{m-1}u)=\sum_{n \geq 0} (q_n \circ p_m)u^{nm}$
and $Q^\perp(-u^{-1}) \cdots Q^\perp(-\omega^{m-1}u^{-1}) = \sum_{n \geq 0} (-1)^n(q_n \circ p_m)^\perp u^{-nm}$. 
Thus the right-hand side of the theorem equals to a constant multiple of the left-hand side. 
Comparing the actions of the constant terms on $1 \in \Gamma$ we conclude that this constant is $1$. 
\hfill $\Box$

\begin{rem}
In the equation (\ref{chounagai}), the product in the middle of the equation is in fact not well-defined. 
However the problem can be avoided by, first calculating the leftmost side 
with replacing the $u$'s appearing on each of $r+1$ factor by $r+1$ distinct variables, obtaining an expression like the rightmost side of (\ref{chounagai}), 
and specializing all variables to $u$. 
\end{rem}

\begin{rem}
As operators $\QB_r$ almost anticommute, it is natural to try constructing, 
in the same way as the construction (\ref{phi2b}) of $\phi_k^{(m)}$, 
an operator $\sum \QB_{a_{m-1}} \cdots \QB_{a_0}$ 
where the sum runs over all $(a_0, \ldots, a_{m-1}) \in \ZZ^m$ with $a_i \equiv i \pmod m$ ($0 \leq i \leq m-1$) and $\sum a_i = (0+\cdots+(m-1))-km$. 
It can be shown, by the same calculation as Theorem \ref{mainthm}, that if $m=2$ this operator commutes with all $p_l^\perp$ ($l: \mathrm{odd}$)
and thus gives a relation of characters. 
However, this operator in fact coincides with $p_{2k-1}^\perp$, 
and the relation obtained thus coincides with ordinary recurrence formula given by expanding $p_l^\perp Q_\lmb$ by $Q$-functions (see eg. \cite[Chap. 10]{HH}). 
\end{rem}

\section{The Brauer algebra and the walled Brauer algebra}
The Brauer algebra was introduced in \cite{Br} in order to describe the centralizer algebra of an orthogonal or a symplectic group
acting on a tensor power of its vector representation. 
The Brauer algebra $D_n(x)$ for a nonnegative integer $n$ and a parameter $x$ (which can be either a complex number or an indeterminate) is defined as follows: 
\begin{itemize}
\item $D_n(x)$ has a basis (over $\CC$ or $\CC(x)$ depending on whether the parameter is a number of an indeterminate)
consisting of all diagrams obtained by connecting $2n$ dots, 
aligned in two rows and $n$ columns, to form $n$ pairs of dots (they are called \textit{Brauer diagrams} or \textit{$n$-diagrams}). 
\item To calculate the product of two Brauer diagrams $a$ and $b$, 
place $a$ on top of $b$ so that the bottom row of $a$ coincides with the top row of $a$, 
remove all closed loops, ignore the dots in the middle row to get another Brauer diagram, 
and multiply it by $x^{\#\text{removed loops}}$. 
\end{itemize}
\begin{center}
\unitlength 0.1in
\begin{picture}( 40.0000,  8.0000)(  2.0000,-10.0000)
%
\special{pn 8}%
\special{pa 200 400}%
\special{pa 400 800}%
\special{fp}%
\special{pa 600 400}%
\special{pa 1000 800}%
\special{fp}%
\special{pa 800 400}%
\special{pa 800 800}%
\special{fp}%
%
\special{pn 8}%
\special{pa 2000 400}%
\special{pa 1800 800}%
\special{fp}%
%
\special{pn 8}%
\special{ar 1400 400 200 100  6.2831853 6.2831853}%
\special{ar 1400 400 200 100  0.0000000 3.1415927}%
%
\special{pn 8}%
\special{ar 1600 400 200 100  6.2831853 6.2831853}%
\special{ar 1600 400 200 100  0.0000000 3.1415927}%
%
\special{pn 8}%
\special{ar 1300 800 100 100  3.1415927 6.2831853}%
%
\special{pn 8}%
\special{ar 1800 800 200 100  3.1415927 6.2831853}%
%
\special{pn 8}%
\special{ar 700 400 300 100  6.2831853 6.2831853}%
\special{ar 700 400 300 100  0.0000000 3.1415927}%
%
\special{pn 8}%
\special{ar 400 800 200 100  3.1415927 6.2831853}%
%
\special{pn 8}%
\special{sh 1}%
\special{ar 1200 800 10 10 0  6.28318530717959E+0000}%
\special{sh 1}%
\special{ar 1400 800 10 10 0  6.28318530717959E+0000}%
\special{sh 1}%
\special{ar 1600 800 10 10 0  6.28318530717959E+0000}%
\special{sh 1}%
\special{ar 1800 800 10 10 0  6.28318530717959E+0000}%
\special{sh 1}%
\special{ar 2000 800 10 10 0  6.28318530717959E+0000}%
\special{sh 1}%
\special{ar 2000 400 10 10 0  6.28318530717959E+0000}%
\special{sh 1}%
\special{ar 1800 400 10 10 0  6.28318530717959E+0000}%
\special{sh 1}%
\special{ar 1600 400 10 10 0  6.28318530717959E+0000}%
\special{sh 1}%
\special{ar 1400 400 10 10 0  6.28318530717959E+0000}%
\special{sh 1}%
\special{ar 1200 400 10 10 0  6.28318530717959E+0000}%
%
\special{pn 8}%
\special{sh 1}%
\special{ar 200 400 10 10 0  6.28318530717959E+0000}%
\special{sh 1}%
\special{ar 400 400 10 10 0  6.28318530717959E+0000}%
\special{sh 1}%
\special{ar 600 400 10 10 0  6.28318530717959E+0000}%
\special{sh 1}%
\special{ar 800 400 10 10 0  6.28318530717959E+0000}%
\special{sh 1}%
\special{ar 1000 400 10 10 0  6.28318530717959E+0000}%
\special{sh 1}%
\special{ar 1000 800 10 10 0  6.28318530717959E+0000}%
\special{sh 1}%
\special{ar 800 800 10 10 0  6.28318530717959E+0000}%
\special{sh 1}%
\special{ar 600 800 10 10 0  6.28318530717959E+0000}%
\special{sh 1}%
\special{ar 400 800 10 10 0  6.28318530717959E+0000}%
\special{sh 1}%
\special{ar 200 800 10 10 0  6.28318530717959E+0000}%
\put(10.7000,-6.5000){\makebox(0,0)[lb]{$\cdot$}}%
\put(20.4000,-6.5000){\makebox(0,0)[lb]{$=$}}%
%
\special{pn 8}%
\special{pa 3000 600}%
\special{pa 2800 1000}%
\special{fp}%
%
\special{pn 8}%
\special{ar 2400 600 200 100  6.2831853 6.2831853}%
\special{ar 2400 600 200 100  0.0000000 3.1415927}%
%
\special{pn 8}%
\special{ar 2600 600 200 100  6.2831853 6.2831853}%
\special{ar 2600 600 200 100  0.0000000 3.1415927}%
%
\special{pn 8}%
\special{ar 2300 1000 100 100  3.1415927 6.2831853}%
%
\special{pn 8}%
\special{ar 2800 1000 200 100  3.1415927 6.2831853}%
%
\special{pn 8}%
\special{sh 1}%
\special{ar 2200 1000 10 10 0  6.28318530717959E+0000}%
\special{sh 1}%
\special{ar 2400 1000 10 10 0  6.28318530717959E+0000}%
\special{sh 1}%
\special{ar 2600 1000 10 10 0  6.28318530717959E+0000}%
\special{sh 1}%
\special{ar 2800 1000 10 10 0  6.28318530717959E+0000}%
\special{sh 1}%
\special{ar 3000 1000 10 10 0  6.28318530717959E+0000}%
\special{sh 1}%
\special{ar 3000 600 10 10 0  6.28318530717959E+0000}%
\special{sh 1}%
\special{ar 2800 600 10 10 0  6.28318530717959E+0000}%
\special{sh 1}%
\special{ar 2600 600 10 10 0  6.28318530717959E+0000}%
\special{sh 1}%
\special{ar 2400 600 10 10 0  6.28318530717959E+0000}%
\special{sh 1}%
\special{ar 2200 600 10 10 0  6.28318530717959E+0000}%
%
\special{pn 8}%
\special{pa 2200 200}%
\special{pa 2400 600}%
\special{fp}%
\special{pa 2600 200}%
\special{pa 3000 600}%
\special{fp}%
\special{pa 2800 200}%
\special{pa 2800 600}%
\special{fp}%
%
\special{pn 8}%
\special{sh 1}%
\special{ar 2200 200 10 10 0  6.28318530717959E+0000}%
\special{sh 1}%
\special{ar 2400 200 10 10 0  6.28318530717959E+0000}%
\special{sh 1}%
\special{ar 2600 200 10 10 0  6.28318530717959E+0000}%
\special{sh 1}%
\special{ar 2800 200 10 10 0  6.28318530717959E+0000}%
\special{sh 1}%
\special{ar 3000 200 10 10 0  6.28318530717959E+0000}%
\special{sh 1}%
\special{ar 3000 600 10 10 0  6.28318530717959E+0000}%
\special{sh 1}%
\special{ar 2800 600 10 10 0  6.28318530717959E+0000}%
\special{sh 1}%
\special{ar 2600 600 10 10 0  6.28318530717959E+0000}%
\special{sh 1}%
\special{ar 2400 600 10 10 0  6.28318530717959E+0000}%
\special{sh 1}%
\special{ar 2200 600 10 10 0  6.28318530717959E+0000}%
%
\special{pn 8}%
\special{ar 2400 600 200 100  3.1415927 6.2831853}%
%
\special{pn 8}%
\special{ar 2700 200 300 100  6.2831853 6.2831853}%
\special{ar 2700 200 300 100  0.0000000 3.1415927}%
\put(30.4000,-6.5000){\makebox(0,0)[lb]{$=$}}%
%
\special{pn 8}%
\special{sh 1}%
\special{ar 3400 400 10 10 0  6.28318530717959E+0000}%
\special{sh 1}%
\special{ar 3600 400 10 10 0  6.28318530717959E+0000}%
\special{sh 1}%
\special{ar 3800 400 10 10 0  6.28318530717959E+0000}%
\special{sh 1}%
\special{ar 4000 400 10 10 0  6.28318530717959E+0000}%
\special{sh 1}%
\special{ar 4200 400 10 10 0  6.28318530717959E+0000}%
\special{sh 1}%
\special{ar 4200 800 10 10 0  6.28318530717959E+0000}%
\special{sh 1}%
\special{ar 4000 800 10 10 0  6.28318530717959E+0000}%
\special{sh 1}%
\special{ar 3800 800 10 10 0  6.28318530717959E+0000}%
\special{sh 1}%
\special{ar 3600 800 10 10 0  6.28318530717959E+0000}%
\special{sh 1}%
\special{ar 3400 800 10 10 0  6.28318530717959E+0000}%
\put(32.0000,-6.5000){\makebox(0,0)[lb]{$x\cdot$}}%
%
\special{pn 8}%
\special{pa 3800 400}%
\special{pa 4000 800}%
\special{fp}%
%
\special{pn 8}%
\special{ar 3700 400 300 100  6.2831853 6.2831853}%
\special{ar 3700 400 300 100  0.0000000 3.1415927}%
%
\special{pn 8}%
\special{ar 3900 400 300 100  6.2831853 6.2831853}%
\special{ar 3900 400 300 100  0.0000000 3.1415927}%
%
\special{pn 8}%
\special{ar 3500 800 100 100  3.1415927 6.2831853}%
%
\special{pn 8}%
\special{ar 4000 800 200 100  3.1415927 6.2831853}%
\end{picture}%

{\scriptsize Figure 5.1: multiplication of two Brauer diagrams}
\end{center}

The Brauer algebra contains the group algebra of $S_n$ as a subalgebra: 
$w \in S_n$ corresponds to the Brauer diagram obtained by connecting the $i$-th dot
on the bottom row with the $w(i)$-th dot on the top row for every $1 \leq i \leq n$. 

The structure and the characters of the Brauer algebra are investigated in \cite{Wz} and \cite{Ram}. 
Here we consider the case where $x$ is an indeterminate. 
In this case, $D_n(x)$ is semisimple, and its irreducible representations are parametrized 
by the partitions $\lmb$ with $|\lmb|=n-2i$, $0 \leq i \leq \lfloor \frac{n}{2} \rfloor$. 
By specializing $x$ to $N=-2m$ (resp. $N=m$) in the $\CC[x]$-form of $D_n(x)$ spanned by the diagram basis, 
one obtains a $\CC$-algebra $D_n(N)$, which surjects onto the centralizers of $Sp(2m,\CC)$ (resp. $O(m,\CC)$) 
acting on the $n$-fold tensor product of the vector representation. 
For each $\lmb$ as above, an irreducible representation of $D_n(-2m)$ can be explicitly constructed for sufficiently large $m$, 
say $H_{2m,n,\lmb}$, as the space of highest vectors of a fixed weight for $Sp(2m,\CC)$ (see \cite{BBL}), 
which allows one to ``read off'' the irreducible representation of $D_n(x)$ labeled by $\lmb$, or of its ``$\CC[x]$-form''
in such a way that setting $x$ to $-2m$ in the representation matrices yields
the irreducible representation $H_{2m,n,\lmb}$. 
Thus the irreducible characters of $D_n(x)$ can be investigated by using these dualities. 
In fact, a character of $D_n(x)$ is determined by its values at certain elements, say, elements of the form $w \otimes e^{\otimes i}$, 
where $w$ is an element in $S_{n-2i}$, 
$e$ is a 2-diagram whose top two dots are connected with each other, and whose bottom two dots are also connected with each other, 
and $\otimes$ denotes the horizontal concatenation of two diagrams. 
Moreover, knowing the character values at permutations is essential: 
character values of $D_n(x)$ at $w \otimes e^{\otimes i}$ can be obtained from character values of $D_{n-2i}(x)$ at $w$ in a simple manner (see remarks at the end of \cite[\S 5]{Ram}). 

Let $\chi^{(n)}_\lmb$ be the irreducible character of $D_n(x)$ corresponding to the partition $\lmb$. 
Then its value $\chi^{(n)}_\lmb(\mu)$ at a permutation with cycle type $\mu$ can be calculated by the following Frobenius-type formula: 
\begin{thm}[{reformulation of \cite[Theorem 6.9]{Ram} in the case of a permutation}]
$\chi^{(n)}_\lmb(\mu)$ equals to the constant term of $A^{\mathrm{br}}_{\mu_1} \cdots A^{\mathrm{br}}_{\mu_l} s_\lmb$, 
where $A^{\mathrm{br}}_{r}=p_r+p_r^\perp+\epsilon_r$, $\epsilon_r=\left\{ \begin{array}{ll} 1&(r: \mathrm{even}) \\ 0&(r: \mathrm{odd}) \end{array} \right.$. 
\end{thm}
By Theorem \ref{bernpm_thm} (and the remark after it), $\phi_{n}^{(m)}$ commutes with $A^{\mathrm{br}}_l$ ($m \nmid l$), 
and, for any $f \in \Lmb$, the constant terms of $A^{\mathrm{br}}_m f$ and $(-\phi_{1}^{(m)}+\epsilon_m)f$ (or $(-\phi_{1}^{(m)}+\epsilon_m\phi_{0}^{(m)})f$)
coincide. 
Thus we we can show, as in the beginning of the proof of Theorem \ref{mainthm}
(with some of the equalities replaced by equalities of constant terms), 
the following: 
\begin{thm}
For $\lmb, \mu$ with $|\lmb|=n-2i$ and $|\mu|=n$, 
if no part of $\mu$ is divisible by $m$ we have
$\chi^{(n)}_\lmb(\mu \cup (m)) = \chi^{(n)}_{\phi(\lmb)}(\mu)$, 
where $\phi=-\phi_{1}^{(m)}+\epsilon_m$ or $-\phi_{1}^{(m)}+\epsilon_m\phi_{0}^{(m)}$. 
\hfill $\Box$
\end{thm}

The walled Brauer algebra $D_{r,s}(x)$, where $r$ and $s$ are nonnegative integers, 
is an algebra that plays the role of the Brauer algebra for the analysis of the representation
of $GL(V)$ on a mixed tensor space $V^{\otimes r}\otimes {V^*}^{\otimes s}$ 
(see \cite{BCHLLS}, \cite[Proof of Lemma 1.2]{Koi}). 
It is the subalgebra of $D_{r+s}(x)$
spanned by all diagrams with the following condition: 
the $1, 2, \ldots, r$-th dots on the top (resp. bottom) row should not be connected with the $(r+1), \ldots, (r+s)$-th dots on the bottom (resp. top) row. 
Such diagrams are often depicted by drawing a ``wall'' between the $r$-th and $(r+1)$-th columns. 

If $x$ is an indeterminate, the $\CC(x)$-algebra $D_{r,s}(x)$ is also semisimple, 
and its its irreducible representations are indexed by the pairs of partitions $[\lmb,\kappa]$  (this bracket has nothing to do with our maya-diagram notation) with $(|\lmb|,|\kappa|)=(r-i,s-i)$, $0 \leq i \leq \min\{r,s\}$. 
Again the irreducible representations of $D_{r,s}(x)$ can be obtained from the irreducible representations of $D_{r,s}(k)$ on the highest weight vectors for $GL(V)$, 
and thus these irreducible characters can be investigated, as in the Brauer-algebra case, using the duality with general linear groups (\cite{BCHLLS, Hal}). 
Again, the whole character table follows from the examination of character values at permutations. 
Let $\chi^{(r,s)}_{[\lmb,\kappa]}$ be the irreducible character of $D_{r,s}(x)$ indexed by $[\lmb, \kappa]$. 
A permutation in $D_{r,s}(x)$ has no lines going across the wall, 
so it makes sense to talk about the ``left side'' (the diagram restricted to $1, \ldots, r$-th columns) and the ``right side'' (the diagram restricted to $(r+1), \ldots, (r+s)$-th columns) of such a permutation. 
Let $\chi^{(r,s)}_{[\lmb,\kappa]}(\mu,\nu)$ be the value of $\chi_{[\lmb,\kappa]}$ at a permutation in $D_{r,s}(x)$ whose left and right sides have cycle types $\mu$ and $\nu$ respectively. 
Then the following Frobenius-type formula holds:
\begin{thm}[{reformulation of the results in \cite[\S 7]{Hal}}]
For $F \in \mathrm{End}(\Lmb)$, denote by $F(X)$ and $F(Y)$ the operators in $\mathrm{End}(\Lmb \otimes \Lmb)$
acting on the first and the second factors of the tensor product respectively. 
Then $\chi^{(r,s)}_{[\lmb,\kappa]}(\mu,\nu)$ equals to the constant term (coefficient of $1 \otimes 1$) of $A^+_{\mu_1} \cdots A^+_{\mu_{l}} A^-_{\nu_1} \cdots A^-_{\nu_{l'}} (s_\lmb \otimes s_\kappa)$, 
where $A^+_i=p_i^\perp(X)+p_i(Y)$ and $A^-_i=p_i(X)+p_i^\perp(Y)$. 
\end{thm}
Again by Theorem \ref{bernpm_thm}, $\phi_n^{(m)}(X)$ and $\phi_n^{(m)}(Y)$ commute with $A^\pm_l$ ($m \nmid l$), 
and, for any $f \in \Lmb \otimes \Lmb$, the constant terms of 
$-\phi_1^{(m)}(X)f$ (resp. $-\phi_1^{(m)}(Y)f$) and $A^+_m$ (resp. $A^-_m$) are equal. 
Thus by the same argument as in the Brauer-algebra case and the symmetric-group case we have: 
\begin{thm}
For $\mu$ and $\nu$ with no parts divisible by $m$ and $\lmb, \kappa$ with $|\mu|=r, |\nu|=s, |\lmb|-|\kappa|=r-s$, 
we have $\chi^{(r,s)}_{[\lmb,\kappa]}(\mu \cup (m),\nu) = \chi_{[-\phi_1^{(m)}(\lmb),\kappa]}(\mu,\nu)$
and $\chi^{(r,s)}_{[\lmb,\kappa]}(\mu,\nu \cup (m)) = \chi_{[\lmb,-\phi_1^{(m)}(\kappa)]}(\mu,\nu)$. 
\hfill $\Box$
\end{thm}

\nocite{*}

\bibliographystyle{plain}
\bibliography{articles_used}

\end{document}